\documentclass[12pt]{article}
\usepackage{amsfonts}
\usepackage{mathrsfs}
\usepackage{amssymb}

\begin{document}
\title{Bergman kernel and complex singularity exponent}
\author{Bo-Yong Chen\footnote{Department of Mathematics, Tongji University, Shanghai 200092, P. R. China (email: boychen{\rm @}mail.tongji.edu.cn)}
 \& Hanjin Lee\footnote{Global Leadership School, Handong Global University, Pohang 791-708,
Republic of Korea
 (email: hlee4jvpc{\rm @}gmail.com)}}
\date{}
\maketitle

\section{ Introduction}

Let $\Omega$ be a bounded $C^\infty$ pseudoconvex domain in
$\mathbb{C}^n$. An important subject in complex analysis is to
understand the boundary behavior of the Bergman kernel $K_\Omega(z)$
of $\Omega$. The special case of strongly pseudoconvex domains is
well-understood through the works of H\"ormander$^{[1]}$,
Diederich$^{[2,\,3]}$ especially the deep work of Fefferman$^{[4]}$.
Precise estimates of $K_\Omega$ for some special pseudoconvex finite
type domains are also available in terms of certain embedded
polydisks$^{[5,\,6]}$. However, it is generally impossible to
construct such embedded polydisks and one only knows
  $K_\Omega\ge
C \delta_\Omega^{-2}$ from the Ohsawa-Takegoshi extension
theorem$^{[7]}$ (see also [8] for a slightly weaker result).
Diederich-Herbort-Ohsawa$^{[9]}$ proved $K_\Omega\ge
C\delta_\Omega^{-2-\epsilon},\,\epsilon>0$ for any pseudoconvex
$\Omega$ of finite type. On the other side, it was Herbort$^{[10]}$
who firstly noticed that the growth of the Bergman kernel may
contain log terms for certain pesudoconvex domain of finite type in
$\mathbb{C}^3$. Recently, a far-reaching generalization of this
phenomenon was made by Kamimoto$^{[11]}$, who in fact obtained a
precise asymptotic expansion of the Bergman kernel non-tangentially
at zero on certain model domain of finite type of form
$\Omega=\{(z,w)\in \mathbb{C}^{n+1}:{\rm Im\,}w>\rho(z)\}$ where
$\rho(e^{i\theta_1}z_1,\ldots,e^{i\theta_n}z_n)=\rho(z)$,
$\theta_j\in \mathbb{R},\,1\le j\le n$, by using the Newton
polyhedron of $\rho$.   In this paper, we shall give a complete
description of the Bergman kernel and metric for the following
important model domains:\vspace{-2mm}
$$
\Omega_F=\{(z,w)\in \mathbb{C}^{n+1}:r(z,w)={\rm
Im\,}w-|F(z)|^2>0\},
$$
where $F=(f_1,\ldots,f_m)$ is a holomorphic map from $\mathbb{C}^n$
to $\mathbb{C}^m$. We always assume $F(0)=0$ for the sake of
simplicity.

To state our results, let us first recall the following

\medskip

\textbf{Definition 1$^{[12]}$.}\ {\it Let $M$ be a complex manifold
and $\rho$ be a measurable function on $M$. For any compact set
$K\subset M$, the \emph{complex singularity exponent} of $\rho$ on
$K$ is defined by}\vspace{-2mm}
$$
c_K(\rho)=\sup\{c\ge 0:|\rho|^{-c} {\ \it is\ } L^2\ {\it on\ a \
neighborhood\ of \ } K\}.\vspace{-2mm}
$$

Roughly speaking, the complex singularity exponent is a holomorphic
invariant which measures the singularity of a function more
precisely than the well-known Lelong number. Given a point $p\in M$,
we write $c_p(\rho)$ instead of $c_{\{p\}}(\rho)$.

\medskip

\textbf{Theorem 1.}\ {\it Let $V\subset\subset U$ be two bounded
Stein neighborhoods of $0\in \mathbb{C}^{n+1}$. Given any
$p_0=(z_0,w_0)\in V\cap
\partial \Omega_F$ and any non-tangential cone $\Lambda\subset \Omega_F$
with vertex at $p_0$, there exists an integer $1\le l\le n$ such
that\vspace{-1mm}
$$
 K_{\Omega_F\cap
U}(p) \asymp  r(p)^{-2-c_{z_0}(|F-F(z_0)|)}|\log r(p)|^{1-l},\ \ \
\forall\,p\in V\cap \Lambda.\vspace{-1mm}
$$
Here $A\asymp B$ means that the ratio $A/B$ is pinched between two
positive constants $($possibly depending on $p_0)$.}

\medskip

\textbf{Remark 1.}\ If in addition $|F(z)|\rightarrow \infty$ as
$|z|\rightarrow \infty$, then the above result holds for
$K_{\Omega_F}$ itself (compared with [13]).

\medskip

To get similar results for the Bergman metric $B_{\Omega_F\cap U}$
and the holomorphic sectional curvature $R_{\Omega_F\cap U}$ of
$B_{\Omega_F\cap U}$, we introduce the following

\medskip

\textbf{Definition 2.}\ {\it Let $\rho,\,\phi$ be non-negative
measurable functions in $\mathbb{C}^n$ and $p\in \mathbb{C}^n$. The
complex singularity exponent with weight $\phi$ of $\rho$ at $p$ is
defined by}
$$
c_p(\rho;\phi)=\sup\{c\ge 0:\phi/\rho^{c} {\ \it is\ } L^2\ {\it on\
a \ neighborhood\ of \ } p\}.
$$

\medskip

\textbf{Theorem 2.}\ {\it Under the conditions of Theorem $1$, there
exist integers $1\le l'_j,l''_j\le n,\ j=1,\ldots,n $ such that
$$
B_{\Omega_F\cap U}(p;X)\asymp
\frac{|X_{n+1}|^2}{r(p)^2}+\sum_{j=1}^n\frac{|X_j|^2}{r(p)^{c_{z_0}(|F-F(z_0)|;|z_j-z_j^0|)-c_{z_0}(|F-F(z_0)|)}|\log
r(p)|^{l'_j-l}},
$$
for all $p\in V\cap \Lambda$ and non-zero vectors $X\in {\mathbb
C}^{n+1}$ and
\begin{eqnarray*}
&& 2-R_{\Omega_F\cap U}[p;(0,\ldots,X_j,\ldots,0)]\\
&&\qquad \asymp
r(p)^{2c_{z_0}(|F-F(z_0)|;|z_j-z_j^0|)-c_{z_0}(|F-F(z_0)|)-c_{z_0}(|F-F(z_0)|;|z_j-z_j^0|^2)}|\log
r(p)|^{2l'_j-l''_j-l},
\end{eqnarray*}
for all $j=1,\ldots,n$. Here the implicit constants are independent
of $X$.}

\medskip

\textbf{Remark 2.}\ Given $c'<c_{z_0} (\rho  )$ and $c''<c_{z_0}
(\rho ; |z_j|^2 )$, we set $c= \frac{c'+c''}{2}$. Then for some
$0<\eta\ll 1$, the Schwarz inequality implies
 $$\int_{B(z_0,\eta)} |\rho|^{-2c} |z_j|^2 dV(z) \le
\bigg( \int_{B(z_0,\eta)} |\rho|^{-2c'}  dV(z) \bigg)^{1/2} \bigg(
\int_{B(z_0, \eta)} |\rho|^{-2c''} |z_j|^4 dV(z) \bigg)^{1/2},
$$
from which it follows that
 $
2c_{z_0}(\rho;|z_j|)\ge c_{z_0}(\rho)+c_{z_0}(\rho;|z_j|^2).
$
Taking $\rho=|F-F(z_0)|$, we see from Theorem 2 that the holomorphic
sectional curvature is either bounded below by a constant or is
asymptotic to $-\infty$ polynomially w.r.t. the Bergman distance in
a non-tangent cone at $p_0$. We conjecture that the latter case
actually can't occur.

\medskip

With the help of some results of Demailly-Koll\'ar$^{[12]}$, we
generalize Theorem 1 to the following

\medskip

\textbf{Theorem 3.}\ {\it Let $\rho\ge 0$ be a log psh function in
${\mathbb C}^n$ with $\rho(0)=0$ and let}
$$
\Omega_\rho:=\{(z,w)\in {\mathbb C}^{n+1}: {\rm Im\,}w>\rho(z) \}.
$$
{\it Let $V\subset\subset U$ be two bounded Stein neighborhoods of
the origin. Then for every $\epsilon>0$, there is a constant
$C_\epsilon\gg 1$ such that }
$$
C_\epsilon^{-1}({\rm Im\,}w)^{-2-2c_0(\rho)+\epsilon} \le
K_{\Omega_\rho\cap U}((0,w))\le C_\epsilon({\rm
Im\,}w)^{-2-2c_0(\rho)-\epsilon}
$$
{\it for $(0,w)\in \Omega_\rho\cap V$}.

\medskip

As an application of Theorem 1, we are able to get asymptotic
estimates for the (Euclidean) volume of sublevel sets with {\it
parameter} $\zeta$
$$
D(U',F,\zeta,r):=\{z\in U':|F(z)-F(\zeta)|<r,\,|\zeta|<r\},\ \ \
U'\subset\subset {\mathbb C}^n
$$
where $F:{\mathbb C}^n\rightarrow {\mathbb C}^m$ is a holomorphic
map such that $F(0)=0$, which might be useful for other purposes.

\medskip

\textbf{Theorem 4.} \ {\it Fix a sufficiently small (Stein)
neighborhood $U'$ of $0\in {\mathbb C}^n$. Then for every
$\epsilon>0$, there is a constant $C_\epsilon=C(U',F,\epsilon)>0$
such that}
$$
C r^{c_0(|F|)}|\log r|^{l-1} \le {\rm Vol\,}(D(U',F,\zeta,r))\le
C_\epsilon r^{c_0(|F|)-\epsilon},\ \ \ {\rm as\ } r\rightarrow 0.
$$
{\it Here $1 \le l\le n$ is an integer coming from the resolution of
singularity of $\{F=0\}$ and $C$ is a constant independent of
$\epsilon$}.

\section{ Preliminaries}\vspace{-3mm}

\subsection{ Bergman invariants and minimum integrals}

 Assume that $\Omega$ is a bounded domain in $\mathbb{C}^n$. Let
$L^2(\Omega)$ denote the space of square-integrable functions on
$\Omega$ and $\|\cdot\|_\Omega$ the corresponding $L^2$-norm. The
Bergman space is given by $H^2(\Omega):=L^2(\Omega)\cap {\cal
O}(\Omega)$ and the Bergman kernel $K_\Omega(z)=\sum_j
|\phi_j(z)|^2$ where $\{\phi_j\}$ is a complete orthogonal basis of
$H^2(\Omega)$. The Bergman metric $B_\Omega(z;X)$ is given by
$\sum_{j,k} g_{j\bar{k}}X_j\bar{X}_k$ where\vspace{-1mm}
$$
g_{j\bar{k}}=\frac{\partial^2 \log K_\Omega(z)}{\partial z_j\partial
\bar{z}_k}.\vspace{-1mm}
$$
The holomorphic sectional curvature of $B_\Omega$ is defined by
$$
R_\Omega(z;X)=B_\Omega(z;X)^{-4}\sum_{h,j,k,l}
R_{\bar{h}jk\bar{l}}\bar{X}_hX_jX_k\bar{X}_l,\ \ \ X\in
\mathbb{C}^n-\{0\},
$$
where
$$
R_{\bar{h}jk\bar{l}}=-\frac{\partial^2 g_{j\bar{h}}}{\partial z_k
\partial \bar{z}_l }+\sum_{\mu,\nu} g^{\nu\bar{\mu}}\frac{\partial
g_{j\bar{\mu}}}{\partial z_k}\frac{\partial g_{\nu\bar{h}}}{\partial
\bar{z}_l},
$$
$g^{\nu\bar{\mu}}$ being the inverse matrix to $g_{j\bar{k}}$. We
define the minimum integrals
\begin{eqnarray*}
&&I^0_\Omega(p)  =  \inf\{\|f\|^2_\Omega: f\in H^2(\Omega),\
f(p)=1\};\\
&&I^1_\Omega(p;X)  =  \inf\bigg\{\|f\|^2_\Omega: f\in H^2(\Omega),\
f(p)=0,\,\sum_j X_j\partial f/\partial z_j (p)=1\bigg\};\\
&&I^2_\Omega(p;X)  =  \inf\bigg\{\|f\|^2_\Omega: f\in H^2(\Omega),\
f(p)=\partial f/\partial z_1(p)=\cdots=\partial f/\partial z_n(p)=
0,\\
&& \qquad\qquad\qquad\quad\sum_{j,k} \partial^2 f/\partial
z_j\partial z_k (p)X_j X_k=1\bigg\}.
\end{eqnarray*}
By the definitions, the minimum integrals increase when the domain
does, and they enjoy the following transformation laws under a
biholomorphic map $\Phi:\Omega_1\rightarrow \Omega_2$:
\begin{eqnarray*}
&&I^0_{\Omega_1}(p)  =  |J_\Phi(p)|^{-2} I^0_{\Omega_2}(\Phi(p)),\\
&&I^{j}_{\Omega_1}(p;X)  =
|J_\Phi(p)|^{-2}I^{j}_{\Omega_2}(\Phi(p);\Phi_\ast X),\ \ \ j=1,2,
\end{eqnarray*}
where $J_\Phi$ denotes the complex Jocobian determinant of $\Phi$.
 The following relationships between the
Bergman invariants and minimum integrals are well-known (see eg.
[14]):
$$
K_\Omega(z)=\frac1{I^0_\Omega(z)},\ \ B_\Omega (z;X)=
\frac{I^0_\Omega(z)}{I^1_\Omega(z;X)},\ \
R_\Omega(z;X)=2-\frac{[I^1_\Omega(z;X)]^2}{I^0_\Omega(z)I^2_\Omega(z;X)}.
$$

\subsection{ Multiplier ideal sheaf}

 One of the most basic concepts on complex analysis is the
multiplier ideal sheaf introduced by Nadel$^{[15]}$.

\medskip

\textbf{Definition 3.}\ {\it Let $\varphi$ be a psh function on a
domain $\Omega$ in $\mathbb{C}^n$. The multiplier ideal sheaf ${\cal
I }(\varphi)\subset {\cal O}_\Omega$ is defined by
$$
\Gamma(U,{\cal I}(\varphi))=\{f\in {\cal
O}_\Omega(U):|f|^2e^{-\varphi}\in L_{\rm loc}^1(U)\}
$$
for every open set $U\subset \Omega$.}\vspace{1mm}
\medskip
For a point $p\in \Omega$ and an integer $l\ge 0$, we set
\begin{eqnarray*}
M_p^l(\Omega) & = & \{\varphi\in PSH(\Omega):\varphi\le 0,\,
\Gamma(U,{\cal I}(\varphi))\subset \mathfrak{M}
_{\Omega,p}^{l+1}\quad {\rm for\ some\ neighborhood\ }U\ {\rm of\
}p\},
\end{eqnarray*}
where $\mathfrak{M}_{\Omega,p}$ denotes the maximal ideal of ${\cal
O}_{\Omega}$ at $p$. We have the following

\medskip

\textbf{Proposition 1.}\ {\it Let $\Omega$ be a pseudoconvex domain
in $\mathbb{C}^n$ and $\varphi\in M_p^l(\Omega)$. Then for any $L^2$
holomorphic function $f$ on $\{\varphi<-1\}$, there is an $L^2$
holomorphic function $\tilde{f}$ on $\Omega$ such that
\begin{eqnarray*}
\frac{\partial^\alpha \tilde{f}}{\partial z_1^{\alpha_1}\cdots
\partial z_n^{\alpha_n} }(p)  =  \frac{\partial^\alpha
f}{\partial z_1^{\alpha_1}\cdots
\partial z_n^{\alpha_n} }(p),\quad
\int_\Omega |\tilde{f}|^2 dV \le  C \int_{\{\varphi<-1\}} |f|^2 dV,
\end{eqnarray*}
for all multi-indexes $\alpha=(\alpha_1,\ldots,\alpha_n)$ with
$|\alpha|\le l$. Here $C$ is a constant depending only on $l$.}

\medskip

\emph{Proof.}\ Take a cut-off function $%
\chi :\mathbb{R}\rightarrow \mathbf{[}0,1\mathbf{]}$ such that $\chi
|_{(-\infty ,-2\log 2]}=1$ and $\chi |_{[-\log 2,+\infty )}=0$. Set
$
\psi=-\log (-\varphi+1).
$
By Donnelly-Fefferman type $L^2$ estimate (see eg. [16]), we can
solve the equation
$
\bar{\partial }u=f\bar{\partial }\chi (\psi)
$
in the weak sense together with the estimate
\begin{eqnarray*}
\int_\Omega|u|^2e^{-\varphi}dV \leq \int_\Omega\left| \bar{\partial
}\chi (\psi)\right| _{\sqrt{-1}\partial \bar{\partial
}\psi}^2|f|^2e^{-\varphi}dV  \leq C_2'\int_{\{\varphi<-1\}} |f|^2 dV
\end{eqnarray*}
because $ \sqrt{-1}\partial \bar{\partial }\psi\geq
\sqrt{-1}\partial\psi \bar{\partial }\psi. $ Here $|\cdot
|_{\sqrt{-1}\partial \bar{\partial }\psi}$ denotes the pointwise
norm with respect to the (singular) metric $\partial \bar{\partial
}\psi$ and $C_2'$ depends only on $n,\,l$ and the choice of $\chi $.
Set
$
\tilde{f}=\chi (\psi) f-u.
$
This $\tilde{f}$ has the desired properties because $u$ is
holomorphic in certain
 neighborhood of $p$.

\subsection{ Pluricomplex Green function}

 Given a domain $\Omega$ in $\mathbb{C}^n$. The pluricomplex Green
function with pole at $p\in \Omega$ is defined by
$$
g_\Omega(z,p)=\sup\{u(z):u\le 0,\,u\in
PSH(\Omega),\,u(z)=\log|z-p|+O(1)\ {\rm near\ }p\}.
$$
One basic fact is that the pluricomplex Green function is decreasing
under holomorphic maps.

\medskip

\textbf{Proposition 2.}\ {\it Let $h$ be a holomorphic map from a
 domain $\Omega\subset \mathbb{C}^n$ to the unit disc $\Delta$.
Then}
$$
\{g_\Omega(\cdot,p)<-1\}\subset \left\{\frac{1-|h(p)|}8\le 1-|h|\le
8(1-|h(p)|)\right\}.
$$
\medskip
\emph{Proof.}\ Observe that
 \begin{eqnarray*}
  -g_\Omega(z,p)
&\leq &-g_\Delta (h(z),h(p))\leq -\log \frac{| h(z)-h(p)| }{|
1-\overline{h(p)}h(z)| }
\\
&= &\frac 12\log \bigg( 1+\frac{( 1-|h(z)|^2) (
1-|h(p)|^2) }{| h(z)-h(p)| ^2}\bigg) \\
&\leq &\frac{( 1-|h(z)|^2) ( 1-|h(p)|^2) }{%
2| h(z)-h(p)| ^2}
 \leq  2\frac{( 1-|h(z)|) ( 1-|h(p)|) }{%
| h(z)-h(p)| ^2}.
\end{eqnarray*}
When $1-|h(z)|\geq 2( 1-|h(p)|)$,
$$
|h(z)-h(p)|\geq 1-|h(z)|- ( 1-|h(p)| ) \geq \frac 12 ( 1-|h(z)| ).
$$
Thus
$$
-g_\Omega(z,p)\leq 8\frac{1-|h(p)|}{1-|h(z)|},
$$
hence
$
 \{g_\Omega(\cdot,p)<-1\}
\subset \{ 1-|h|\leq 8( 1-|h(p)|)\}. $ A similar argument implies
$$
  \{g_\Omega(\cdot,p)<-1\}
\subset \left\{  1-|h|\geq \frac{1-|h(p)|}8 \right\} .
$$

\subsection{ Resolution of singularity}

The following fundamental theorem  will play an essential role in
our proofs.

\medskip

\textbf{Hironaka's Theorem$^{[17]}$.}\ {\it Let $M$ be a complex
manifold and ${\cal I}\subset {\cal O}_M$ be a coherent ideal sheaf.
Then there is a log canonical resolution of ${\cal I}$, i.e., there
exists a proper bimeromorphic morphism $\mu$ from a complex manifold
$\widetilde{M}$ to $M$ such that $\mu^\ast{\cal I}={\cal
O}_{\widetilde{M}}(-D)$ is an invertible sheaf associated to a
divisor $D$ with $D+E$ being simple normal crossings. Here $E$
denotes the exceptional divisor of $\mu$.}\vspace{1mm}

\medskip

Let $F:\mathbb{C}^n\rightarrow \mathbb{C}^m$ be a holomorphic map
with $F(0)=0$, $\mu:\widetilde{M}\rightarrow \mathbb{C}^n$ a log
canonical resolution of ${\cal I}=F^{-1}(0)$. Let $U$ be a
sufficiently small open neighborhood of $0$. By the Jacobian formula
for a change of variable, we have
$$
\int_U
|F(\zeta)|^{-2c}dV_\zeta=\int_{\mu^{-1}(U)}|F\circ\mu(z)|^{-2c}|J_\mu(z)|^2d\widetilde{V}_z,
$$
where $dV_\zeta, \,d\widetilde{V}_z$ are volume elements of
$\mathbb{C}^n,\,\widetilde{M}$ respectively. According to Hironaka's
theorem, this integral is given by a finite number of integrals of
form
$$
\int_{\mu^{-1}(U)\cap
\widetilde{U}}|F\circ\mu(z)|^{-2c}|J_\mu(z)|^2d\widetilde{V}_z
$$
over suitable coordinate charts $\widetilde{U}\subset \widetilde{M}$
on which one has
$$
|F\circ\mu(z)|\asymp \prod_{j=1}^n |z_j|^{a_j},\ \ \ \
|J_\mu(z)|\asymp \prod_{j=1}^n |z_j|^{b_j},
$$
for certain  non-negative integers $a_j,\,b_j$. Since
$|F\circ\mu|^{-2c}|J_\mu|^2$ is $L^1$ on $\mu^{-1}(U)\cap
\widetilde{U}$ if and only if $c<(b_j+1)/a_j$ for all $1\le j\le n$,
it follows that $ c_0(|F|)=\min \{(b_j+1)/a_j\}, $ where the minimum
is taken over all $a_j,b_j$ associated to those
  $\widetilde{U}$. In particular, $c_0(|F|)$ is a
rational number.
\medskip
 Analogously, one can compute $c_0(|F|;|\zeta_k|^{\tau})$. We assume
 $k=1$ and write
 $$
 \mu=(\mu_1,\ldots,\mu_n),\ \ \ |\mu_1(z)|\asymp \prod_{j=1}^n
 |z_j|^{c_j},
 $$
 for certain non-negative integers $c_j$ in above $\widetilde{U}$. Then
 $$
\int_U \frac{|\zeta_1|^{2\tau}}{|F(\zeta)|^{2c}}dV_\zeta
=\int_{\mu^{-1}(U)}\frac{|\mu_1(z)|^{2\tau}|J_\mu(z)|^2}{|F\circ\mu(z)|^{2c}}d\widetilde{V}_z,
$$
 therefore,
 $
 c_0(|F|;|\zeta_1|^\tau)=\min \{(b_j+c_j\tau+1)/a_j\}.
 $

\section{ Proof of the theorems}

Without loss of generality, we assume $U=U'\times U''$ with Stein
open sets $U'\subset \mathbb{C}^n$, and $U''\subset \mathbb{C}$.
Given $w\in \mathbb{C}$, we set
\begin{eqnarray*}
\Omega^1_w &=& \{z\in U':|F(z)|^2< 9\,{\rm Im\,}w\},\quad \Omega^2_w
 =  \bigg\{z\in U':|F(z)|^2< \frac19\,{\rm
Im\,}w\bigg\},\\
S_w & = & \bigg\{\zeta\in \mathbb{C}: \frac19\,{\rm Im\,}w<{\rm
Im\,}\zeta<9\,{\rm Im\,}w \bigg\}.
\end{eqnarray*}

\medskip

\textbf{Lemma 1.}\ {\it There exists a universal constant $C>0$ such
that for any $p=(z,w)\in (\Omega^2_{w}\times S_{w})\cap U,$}
\begin{eqnarray}
&&I^0_{\Omega^2_w\times S_w}(p)  \le  I^0_{\Omega_F\cap U}(p)
\le  C\,I^0_{\Omega^1_w\times S_w}(p),\\
&&I^j_{\Omega^2_w\times S_w}(p;X)  \le  I^j_{\Omega_F\cap U}(p;X)
\le  C\,I^j_{\Omega^1_w\times S_w}(p;X),\ \ \ j=1,2.
\end{eqnarray}

\medskip

\emph{Proof.}\ The first inequalities in (1) and (2) follow directly
from the definitions of the minimum integrals.  Choosing
$h(z,\zeta)=e^{i\epsilon \zeta/{\rm Im\,}w}$ for sufficiently small
$\epsilon$ and applying Proposition 2, we obtain
$$
\{g_{\Omega_F\cap U}(\cdot,p)<-1\}\subset \Omega^1_w\times S_w.
$$
Since $(2n+l)g_{\Omega_F\cap U}(\cdot,p)\in M_p^l(\Omega_F\cap U)$,
the second inequalities in (1) and (2) follow from Proposition 1 and
the definitions of the minimum integrals.
\medskip
Note that $\Phi(\zeta)=e^{i\zeta/{\rm Im\,}w}$ maps the strip $S_w$
biholomorphically to the ring $\{e^{-9}<|\zeta|<e^{-1/9}\}$ with
$|\Phi(w)|=e^{-1}$, thus
$$
I^0_{S_w}(w)  \asymp  ({\rm Im\,}w)^2,\ \  I^1_{S_w}(w;X_{n+1})
\asymp  \frac{({\rm Im\,}w)^2}{|X_{n+1}|^2},\ \ I^2_{S_w}(w;X_{n+1})
\asymp  \frac{({\rm Im\,}w)^2}{|X_{n+1}|^4},
$$
consequently,
\begin{eqnarray}
&& I^0_{\Omega^k_w\times S_w}(p)\asymp ({\rm Im\,}w)^2
I^0_{\Omega^k_w}(z),\quad I^j_{\Omega^{k}_w\times
S_w}[p;(X',0)]\asymp ({\rm Im\,}w)^2I^j_{\Omega^k_w}(z;X'),\\
&& I^1_{\Omega^{k}_w\times S_w}[p;(0,X_{n+1})]\asymp \frac{({\rm
Im\,}w)^2}{|X_{n+1}|^2}I^0_{\Omega^k_w}(z),
\end{eqnarray}
for $j,k=1,2$.

\medskip

\textbf{Remark 3.}\ Actually, the above conclusions still hold if
one replaces $|F|^2$ in $\Omega_F$ by any non-negative psh function.

\medskip

\textbf{Lemma 2.}\ {\it For $\tau=0,1,2$, there exist integers $1\le
m_\tau \le n$ such that
\begin{equation}
\int_{\Omega^j_w} |\zeta_1|^{2\tau} dV_\zeta  \asymp  ({\rm
Im\,}w)^{c_0(|F|;|\zeta_1|^\tau)}|\log ({\rm Im\,}w)|^{m_\tau-1},\ \
j=1,2.
\end{equation}
Moreover, $m_\tau$ can be computed through a log canonical
resolution}.

\medskip

\emph{Proof.}\ Let $\mu:\widetilde{M}\rightarrow \mathbb{C}^n$ be a
log canonical resolution of $F^{-1}(0)$. Fix a small neighborhood
$U$ of $0\in \mathbb{C}^n$. Similar to Subsection 2.4, the integral
$ \int_{\{\zeta\in U:|F(\zeta)|<t\}}|\zeta_1|^{2\tau} dV_\zeta $
 is then
determined by integrals of form
$$
\int_{\mu^{-1}(U)\cap \{z\in \widetilde{U}:\prod_{j=1}^n
|z_j|^{a_j}<t\}}\prod_{j=1}^n |z_j|^{2b_j+2c_j\tau}d\widetilde{V}_z
$$
over finite suitable coordinate charts $\widetilde{U}\subset
\widetilde{M}$. Using polar coordinates,  it suffices to estimate
the integrals
$$
\tilde{I}(t)=\int_{\{y\in {\rm I}^n\cap
\mathbb{R}^n_+:y_1^{a_1}\cdots y_n^{a_n}<t\}}
y_1^{2b_1+2c_1\tau+1}\cdots y_n^{2b_n+2c_n\tau+1}dy_1\wedge\cdots
\wedge dy_n,
$$
where ${\it I}=[0,1]$. Set $h(y)=y_1^{a_1}\cdots y_n^{a_n}$. We
consider the following two useful integrals:
\begin{eqnarray*}
&&H(s)  =  \int_{{\rm I}^n} h(y)^{-s} y_1^{2b_1+2c_1\tau+1}\cdots
y_n^{2b_n+2c_n\tau+1}dy_1\wedge\cdots \wedge dy_n,\\
&&G(t)   =  \int_{\{y\in {\rm I}^n:h(y)=t\}}
y_1^{2b_1+2c_1\tau+1}\cdots y_n^{2b_n+2c_n\tau+1}dy_1\wedge \cdots
\wedge dy_n/dh.
\end{eqnarray*}
A direct computation shows
\begin{eqnarray*}
H(s)   =   C_s \prod_{j=1}^n[(b_j+c_j\tau+1)/a_j-s]^{-1}
          =   C_s \prod_{j\in J}[(b_j+c_j\tau+1)/a_j-s]^{-l_j}
\end{eqnarray*}
provided $s<(b_j+c_j\tau+1)/a_j,\ \forall\,1\le j\le n$. Here
$C_s\in (0,\infty)$ and $l_j$ is the multiplicity of $H(s)$ at the
pole $s=(b_j+c_j\tau+1)/a_j$. As
$$
dy_1\wedge\cdots\wedge dy_n=(a_1h)^{-1}y_1 dh\wedge
dy_2\wedge\cdots\wedge dy_n,
$$
it follows that
$$
G(t)=\frac1{a_1}t^{\frac1{a_1}-1}\int_{{\rm I}^{n-1}\cap
\mathbb{R}^{n-1}_+\cap\{ y_2^{a_2}\cdots y_n^{a_n}\ge
t\}}y_2^{2b_2+2c_2\tau+1-\frac{a_2}{a_1}}\cdots
y_n^{2b_n+2c_n\tau+1-\frac{a_n}{a_1}}dy_2\wedge\cdots\wedge dy_n.
$$
Therefore, $G(t)$ must be of form
$$
\sum_{j=1}^{j_0}\sum_{k=1}^{k_j(\le n)} c_{j,k} t^{\nu_j-1}(-\log
t)^{k-1}
$$
for sufficiently small $t$. Noting that
$
H(s)=\int_0^\infty t^{-s} G(t) dt,
$
by comparing the multiplicities of poles of both sides, we obtain
$$
G(t)\asymp t^{\beta-1}(-\log t)^{\alpha-1},\ \ \ t\rightarrow 0_+,
$$
where
$$
\beta=\min_{1\le j\le n} \{(b_j+c_j\tau+1)/a_j\},\ \ \ \
\alpha=\max_{(b_j+c_j\tau+1)/a_j=\beta,j\in J} l_j.
$$
It follows that
$$
\tilde{I}(t)=\int_0^t G(t')dt'\asymp t^{\beta}(-\log t)^{\alpha-1},\
\ \ t\rightarrow 0_+.
$$
Since by Subsection 2.4,
$
c_0(|F|;|\zeta_1|^\tau)=\min\{\beta\},
$
the lemma is verified with $ l=\max\{\alpha\}, $ where the minimum
is taken over all $\widetilde{U}$ and the maximum is taken over
those $\widetilde{U}$ such that the associated $\beta$ is equal to
$c_0(|F|;|\zeta_1|^\tau)$.

\medskip

\textbf{Remark 4.}\ By the definitions, the minimum integrals
satisfy
\begin{eqnarray}
&&I^0_{\Omega^1_w}(z)\le \int_{\Omega^1_w} dV_\zeta,\ \
I^1_{\Omega^1_w}(z;(X_1,0,\ldots,0))\le |X_1|^{-2}\int_{\Omega^1_w}
|\zeta_1-z_1|^{2} dV_\zeta,
\\&&I^2_{\Omega^1_w}(z;(X_1,0,\ldots,0))\le
|X_1|^{-4}\int_{\Omega^1_w} |\zeta_1-z_1|^{4} dV_\zeta.
\end{eqnarray}

\medskip

\textbf{Lemma 3.}\ {\it Let $m_\tau,\,\tau=0,1,2$ be as in Lemma
$2$. Then}
\begin{eqnarray}
&&I^0_{\Omega^2_w}(0)\ge C\int_{\Omega^2_w} dV_\zeta,\ \
I^1_{\Omega^2_w}(0;X_1)\ge C |X_1|^{-2}\int_{\Omega^2_w}
|\zeta_1|^{2} dV_\zeta, \\&&I^2_{\Omega^2_w}(0;X_1)\ge
C|X_1|^{-4}\int_{\Omega^2_w} |\zeta_1|^{4} dV_\zeta.
\end{eqnarray}
\medskip
\emph{Proof.}\ Let $\mu$ be as above. For any $f\in H^2(\Omega^2_w)$
with $f(0)=1$, we have
$$
\int_{\Omega^2_w}
|f(\zeta)|^{2}dV_\zeta=\int_{\mu^{-1}(\Omega^2_w)}|f\circ\mu(z)|^2|J_\mu(z)|^2d\widetilde{V}_z,
$$
which is then given by integrals of form\vspace{-2mm}
$$
\int_{\mu^{-1}(U)\cap \{z\in \widetilde{U}:\prod_{j=1}^n
|z_j|^{a_j}<t\}} |f\circ\mu(z)|^2\prod_{j=1}^n
|z_j|^{2b_j}d\widetilde{V}_z\vspace{-2mm}
$$
over finite coordinate charts $\widetilde{U}\subset \widetilde{M}$,
where $t={\rm Im\,}w$. Without loss of generality, we assume that
$\widetilde{U}$ contains the unit polydisc $\Delta^n$. As $f\circ
\mu$ is holomorphic on the Reinhardt domain $\{z\in
\Delta^n:\prod_{j=1}^n |z_j|^{a_j}<t\}$, it has an
expansion\vspace{-2mm}
$$
f\circ\mu(z)=1+\sum_{\gamma\neq 0} c_\gamma z^\gamma,\ \ \
\gamma=(\gamma_1,\ldots,\gamma_n),\vspace{-2mm}
$$
because $f\circ \mu (0)=1$.  Therefore,\vspace{-2mm}
\begin{eqnarray*}
&& \int_{\{z\in \Delta^n:\prod_{j=1}^n |z_j|^{a_j}<t\}}
|f\circ\mu(z)|^2\prod_{j=1}^n |z_j|^{2b_j}d\widetilde{V}_z \\
& &\quad=  \int_{\{z\in \Delta^n:\prod_{j=1}^n |z_j|^{a_j}<t\}}
\prod_{j=1}^n |z_j|^{2b_j}d\widetilde{V}_z+\sum_{\gamma\neq 0}
|c_\gamma|^2\int_{\{z\in \Delta^n:\prod_{j=1}^n |z_j|^{a_j}<t\}}
\prod_{j=1}^n |z_j|^{2b_j+2\gamma_j}d\widetilde{V}_z,
\end{eqnarray*}
consequently, $\|f\|^2_{\Omega^2_w}$ dominates all integrals of form
$$
\int_{\{z\in \Delta^n:\prod_{j=1}^n |z_j|^{a_j}<t\}} \prod_{j=1}^n
|z_j|^{2b_j}d\widetilde{V}_z,
$$
hence the integral $\int_{\Omega^2_w} dV_\zeta$. Next, let $g$ be a
given $L^2$ holomorphic function on $\Omega^2_w$ satisfying $g(0)=0$
and $X_1\partial g/\partial \zeta_1(0)=1$. Then
$$
g\circ\mu(z)=\sum_{\gamma\neq 0} \delta_\gamma z^\gamma,\ \ \ {\rm
for\ } z\in \Delta^n,\quad\prod_{j=1}^n |z_j|^{a_j}<t,\vspace{-1mm}
$$
because $g\circ\mu(0)=0$. Since $\mu$ is locally biholomorphic on
$\widetilde{M}-D-E$ where $\mu^\ast\circ F^{-1}(0)={\cal
O}_{\tilde{M}}(-D)$ and $E$ is the exceptional divisor of $\mu$, it
follows that for any $k\ge 2$, $|\mu_k(z)|\asymp
|z_1|^{\gamma_1}\cdots |z_n|^{\gamma_n}$ with
$(\gamma_1,\ldots,\gamma_n)\neq (c_1,\ldots,c_n)$, consequently,
$$
\left.\frac{\partial^{|c_1|+\cdots+|c_n|}\mu_k}{\partial
z_1^{c_1}\cdots\partial z_n^{c_n}}\right|_{z=0}=0,\quad k\ge 2,
$$
which implies
$$
 \delta_{c_1\cdots
c_n}=\left.\frac{\partial^{|c_1|+\cdots+|c_n|}g\circ\mu}{\partial
z_1^{c_1}\cdots\partial z_n^{c_n}}\right|_{z=0}=\frac{\partial
g}{\partial \zeta_1}(0)=\frac1{X_1}.
$$
Thus $\|g\|^2_{\Omega^2_w}$ must dominate the $|X_1|^{-2}$ multiple
of the sum of integrals\vspace{-2mm}
$$
\int_{\{z\in \Delta^n:\prod_{j=1}^n |z_j|^{a_j}<t\}} \prod_{j=1}^n
|z_j|^{2b_j+2c_j}d\widetilde{V}_z\vspace{-1mm}
$$
over those coordinate charts $\widetilde{U}\subset \widetilde{M}$,
which is equivalent to $\int_{\Omega^2_w} |\zeta_1|^{2} dV_\zeta$.
Inequality (9) can be verified similarly.

\medskip

\emph{Proof of Theorems $1,2$.}\ First we assume $p_0=0$. The
theorems follow directly from Subsection 2.1 and equations (1)--(9).
The general case follows from the transformation laws of the minimum
integrals by noting that the transformation $\Phi$ maps $\Omega_F$
biholomorphically to the domain
$$
\Omega_F'=\{(z',w')\in \mathbb{C}^{n+1}:{\rm
Im\,}w'>|F(z_0+z')-F(z_0)|^2\},
$$
where
\begin{eqnarray}
z'  =  z-z_0,\quad w'  =  w-w_0-2i\sum_{k=1}^m
\overline{f_k(z_0)}(f_k(z)-f_k(z_0)).
\end{eqnarray}

\section{ Newton polyhedron and holomorphic sectional curvature}

There is an effective way to compute the complex singularity
exponent in terms of the Newton polyhedron by using toric resolution
of singularity. Given a power series $f(z)=\sum_\alpha c_\gamma
z^\gamma$ over $\mathbb{C}^n$ with $f(0)=0$. The Newton polyhedron
$\Gamma(f)$ of $f$ is the convex hull of the set $\bigcup
(\gamma+\mathbb{R}^n_+)$, where $\mathbb{R}^n_+$ is the positive
octant and the union is taken over multi-indexes $\gamma$ such that
$c_{\gamma}\neq 0$. We associate every compact face $\Delta$ of
$\Gamma(f)$ with a polynomial $f_\Delta (z)=\sum_{\gamma\in \Delta}
c_\gamma z^\gamma$. We say that $f$ is non-degenerate on $\Delta$ if
$d f_\Delta=0$ has no solution in $(\mathbb{C}^\ast)^n$ where
$\mathbb{C}^\ast=\mathbb{C}-\{0\}$. We say that $f$ is
non-degenerate if $f_\Delta$ is non-degenerate for any $\Delta$.
Assume that the line $\{((1+\tau)t,t,\ldots,t)\}$ intersects the
boundary of $\Gamma(f)$ at point
$Q_\tau=((1+\tau)d_\tau,d_\tau,\ldots,d_\tau)$. Let $\hat{m}_\tau$
denote the number of compact faces of $\Gamma(f)$ containing the
point $Q_\tau$. Then the following fact is well known (see eg.
[18]).

\medskip

\textbf{Proposition 3.}\ Assume $f$ is a non-degenerate entire
function in $\mathbb{C}^n$ with $d_0>1$. Then
\begin{eqnarray}
&&c_0(|f|)=1/d_0,\ \ \ c_0(|f|,|z_1|^\tau)=1/d_\tau,\ \ \tau=1,2;\\
&&m_\tau=\min\{\hat{m}_\tau,n\},\ \ \tau=0,1,2.
\end{eqnarray}

\medskip

As Subsection 2.1 shows, the holomorphic sectional curvature is
always bounded above by 2. It is natural to ask whether this upper
bound is optimal. It is also interesting to ask whether the
holomorphic sectional curvature is bounded below. For the special
case $\Omega_f=\{(z,w)\in \mathbb{C}^{n+1}:{\rm Im\,}w>|f(z)|^2\}$
with $f$ being a non-degenerate entire function, we have the
following self-contained characterization.

\medskip

\textbf{Proposition 4.}\ {\it Let $f$ be as in Proposition $3$,
 $U$ be a
bounded Stein neighborhood of $0\in \mathbb{C}^{n+1}$.

{\rm (i)} If $Q_\tau,\,\tau=0,1,2$ are not contained in the same
supporting hyperplane of $\Gamma(f)$, then $R_{\Omega_f\cap
U}(p;(X_1,0,\ldots,0))\rightarrow 2 $ as $p\rightarrow 0$
non-tangentially.

{\rm (ii)}~Otherwise, $ R_{\Omega_f\cap U}(p;(X_1,0,\ldots,0)) $ is
bounded below by a constant in a non-tangent cone $\Lambda$ at $0$.}

\medskip

\textbf{Example.}\ Take $f(z_1,z_2)=z_1^4+z_1^2z_2+z_1z_2^2+z_2^4$.
It is easy to verify that $Q_\tau,\,\tau=0,1,2$ are not contained in
any hyperplane supporting the Newton polyhedron.
\medskip
\emph{Proof.}\ Suppose that $H$ is a hyperplane which supports
$\Gamma(f)$ at $Q_1$, given by the equation $ \sum_{k=1}^n
x_k/a_k=1$ where $a_k>0$. Since $\Gamma(f)$ is convex,
$$
\frac1 {d_0} \le \sum_{k=1}^n \frac1 {a_k},\quad  \frac1 {d_2} \le
\sum_{k=2}^n \frac1 {a_k}+\frac3 {a_1},
$$
hence
$$
\frac1 {d_2}+\frac1 {d_0}\le 2\bigg(\sum_{k=2}^n \frac1 {a_k}+\frac2
{a_1}\bigg)=\frac2 {d_1}.
$$
It follows that the equality holds if and only if
$Q_\tau,\,\tau=0,1,2$ are all contained in $H$. Combining Theorem 2
with Proposition 3, (i) is verified. On the other hand, if
$Q_\tau,\,\tau=0,1,2$ are all contained in a supporting hyperplane
supporting $\Gamma(f)$, then any compact face of $\Gamma(f)$
containing $Q_1$ must contain $Q_0,\,Q_2$, consequently, $2m_1\le
m_2+m_0$ and (ii) follows from Theorem 2.

\section{ Proofs of Theorems 3,4}

\emph{Proof of Theorem 3.}\ We keep the notions as above. Without
loss of generality, we assume $U=U'\times U''$ with Stein open sets
$U'\subset {\mathbb C}^n$, and $U''\subset {\mathbb C}$. Given $w\in
U''$, we set
\begin{eqnarray*}
\Omega^1_w &=& \{z\in U':\rho(z)< 9\,{\rm Im\,}w\}\\
\Omega^2_w &=& \left\{z\in U':\rho(z)< \frac19\,{\rm
Im\,}w\right\}\\
S_w & = & \left\{\zeta\in {\mathbb C}: \frac19\,{\rm Im\,}w<{\rm
Im\,}\zeta<9\,{\rm Im\,}w \right\}.
\end{eqnarray*}
Without any change of the above argument, we can prove the following
$$
C K_{\Omega^1_w\times S_w}((0,w))\le K_{\Omega_\rho\cap U}((0,w))\le
K_{\Omega^2_w\times S_w}((0,w)),\ \ \ w\in U''
$$
for suitable constant $C>0$. It is easy to see
\begin{equation}
K_{\Omega^j_w\times S_w}((0,w))\asymp ({\rm
Im\,}w)^{-2}K_{\Omega^j_w}(0),\ \ \ j=1,2.
\end{equation}
Hence it suffices to estimate $K_{D_r}(0)$ in terms of certain power
of $r$ where
$$
D_r:= \{z\in U':\rho(z)< r\},\ \ \ r\ll 1.
$$
First we have the trivial inequality
\begin{equation}
K_{D_r}(0)\ge \frac1{{\rm Vol\,}(D_r)}.
\end{equation}
Since $\varphi:=\log \rho$ is psh, we infer from Proposition 4.3 (1)
in [12] that for all positive real number $c<c_0(\rho)$ there is an
estimate
\begin{equation}
{\rm Vol\,}(D_r)\le C(c) r^{2c}.
\end{equation}
(Here we remark that the authors of [12] use the notion
$c_0(\varphi)$ for $c_0(e^\varphi)$ when $\varphi$ is psh). The
first inequality in Theorem 3 is then an consequence of (13)--(15).
For the second inequality, we shall use the celebrated Demailly's
approximation: there is a constant $C>0$ independent of $m$ and
$\varphi$ such that
$$
\varphi(z)-\frac{C}m \le \psi_m(z):=\frac1{2m} \log \sum_k
|g_{m,k}(z)|^2
$$
where $\{g_{m,k}\}$ is an orthonormal basis of ${\cal
H}_{m\varphi}(U')$, the Hilbert space of holomorphic functions $f$
on $U'$ such that
$$
\int_{U'}|f|^2 e^{-2m\varphi}dV<\infty.
$$
It follows from the strong Noetherian property that there exists an
integer $k_0(m)$ and a constant $C_{m,1}>0$ such that
$$
\psi_m-C_{m,1} \le \psi_{m,0}:= \frac1{2m} \log \sum_{0\le k\le
k_0(m)}|g_{m,k}|^2\le \psi_m \ \ \ {\rm on\ } U'.
$$
Thus
$$
D_r=\{z\in U':\varphi(z)<\log r\}\supset \{z\in
U':\psi_{m,0}(z)<\log r -C_{m,2}\}:=D_{r,m}
$$
for some constant $C_{m,2}>0$. Note that
\begin{equation}
K_{D_r}(0)\le K_{D_{r,m}}(0).
\end{equation}
By Theorem 1,
\begin{equation}
K_{D_{r,m}}(0)\le C_{m,3} r^{-2c_0(e^{\psi_{m,0}})}|\log r|^{1-l_m}
\end{equation}
where $1\le l_m\le n$ is certain integer coming from the resolution
of the singularity of $\psi_{m,0}$. As
$$
c_0(e^{\psi_{m,0}})=c_0(e^{\psi_{m}})\rightarrow
c_0(e^{\varphi})=c_0(\rho)
$$
by Theorem 4.2 (3) in [12], the second inequality  follows from
(13), (16) and (17). The proof is complete.

\medskip

\emph{Proof of Theorem 4.}\ We fix a sufficiently small Stein
neighborhood $U=U'\times U''$ of $0\in {\mathbb C}^{n+1}$. Fix
arbitrary $(z_0,w_0)\in \Omega_F\cap U$ such that ${\rm
Im\,}w_0-|F(z_0)|^2=r/9$ and $|z_0|\le r$. Take a holomorphic
transformation  $\Phi$ as (10), we have
$$
K_{\Omega_F\cap U}((z_0,w_0))  =   K_{\Omega_F'\cap
\Phi(U)}(\Phi(z_0,w_0))  \ge  C\frac1{{\rm Vol\,}(D(U',F,z_0,r))}
$$
where the second inequality follows from Lemma 1  with $C$ a
universal constant.

Now $(z_0,w_0)$ lies in a non-tangential cone with vertex at the
origin, we have
$$
K_{\Omega_F\cap U}((z_0,w_0))\asymp r^{-c_0(|F|)}|\log r|^{1-l}
$$
by Theorem 1. Thus we get the first inequality in Theorem 4.

For the second inequality, we use Lemma 3.2 (2) in [3] that for any
$c<c_0(|F|)$, there exists a neighborhood $U'_c$ of $0$ such that
$$
\int_{U'} |F(z)-F(z_0)|^{-c}dV(z)\le C(c), \ \ \ z_0\in U'_c
$$
(Shrinking $U'$ if necessary). Since
$$
\int_{U'} |F(z)-F(z_0)|^{-c}dV(z)\ge \frac1{r^c}{\rm
Vol\,}(D(U',F,z_0,r))
$$
for $r\ll 1$, we are done.

\section{ Remarks and questions}

\textbf{Remark 5.}\ Generally, the conclusion of Theorem 1 fails for
domains $\Omega_\rho=\{(z,w)\in \mathbb{C}^{n+1}: {\rm
Im\,}w>\rho(z)\}$ when $\rho$ is a non-negative real-analytic psh
function in $\mathbb{C}^n$. Consider a power series $\rho=\sum
c_\gamma x^{2\gamma_1}y^{2\gamma_2}$ in $\mathbb{R}^2$ with all
$c_\gamma\ge 0$. Write $z=x+iy$. Then $\rho$ is a non-negative
subharmonic function. Assume that the Newton polyhedron
$\Gamma(\rho)$ of $\rho$ (over $\mathbb{R}$) intersects the $x,y$
axes. Then $c_0(\rho)=1/d_0$ where $d_0$ denotes the distance to
$\Gamma(\rho)$. Set $\delta=\inf \{|\gamma|:c_\gamma>0\}$. Note that
the domains
$$
D_t=\{z\in \mathbb{C}:\rho(z)<t\},\ \ \ D'_t=\{\zeta=(\xi,\eta)\in
\mathbb{C}:\rho(t^{\frac1{2\delta}}\zeta)<t\}
$$
are biholomorphically equivalent, therefore $
K_{D_t}(0)=t^{-\frac1{\delta}}K_{D'_t}(0)\asymp t^{-\frac1{\delta}},
$ because $D_t'$ is pinched between two planar domains
$$
\bigg\{\zeta\in \mathbb{C}:
\sum_{c_\gamma>0,|\gamma|=\delta}\xi^{2\gamma_1}\eta^{2\gamma_2}<\epsilon,|\zeta|<\epsilon\bigg\}
,\ \ \bigg\{\zeta\in \mathbb{C}:
\sum_{c_\gamma>0,|\gamma|=\delta}\xi^{2\gamma_1}\eta^{2\gamma_2}<\frac1
\epsilon\bigg\}
$$
for some $0<\epsilon\ll 1$ independent of $t$. By the remarks under
Theorem 1 and Lemma 1, we have
$$
K_{\Omega_\rho}((0,w))\asymp ({\rm Im\,}w)^{-2-\frac1{\delta}}.
$$
Nevertheless, $c_0(\rho)\neq 1/\delta$ (i.e., $d_0\neq \delta$) in
general, for instance, one can take
$\rho(z)=x^8+x^4y^2+x^2y^6+y^{10}$, then $d_0=10/3>3= \delta$.

\medskip

\textbf{Remark 6.}\ Given a point $p_0\in \mathbb{C}^n$. Let ${\cal
S}$ denote the space of all bounded $C^2$ pseudoconvex domains in
$\mathbb{C}^n$ whose boundary contains $p_0$. For any $\Omega\in
{\cal S}$, we define the growth exponent of the Bergman kernel of
$\Omega$ at $p_0$ by
$$
b_{p_0}(\Omega)=\sup\Big\{b\ge 0: \lim_{p\rightarrow p_0,p\in
\Lambda}\delta_\Omega(p)^bK_\Omega(p)=\infty\Big\},
$$
where $\Lambda$ is some non-tangent cone at $p_0$ and
$\delta_\Omega$ denotes the boundary distance function. Clearly,
$b_{p_0}$ defines a map from ${\cal S}$ to $[2,n+1]$. Note that for
those domains considered in Theorem 1, the values of $b_{p_0}$ are
always rational numbers. Thus it is natural to ask

\medskip

\textbf{Question 1.}\ Is the image of $b_{p_0}$ dense in $[2,n+1]$?
Is $b_{p_0}$ surjective?

\medskip

\textbf{Remark 7.}\ We can't get global uniform estimates of the
Bergman invariants as in the case of strongly pseudoconvex domains
or finite type domains in $\mathbb{C}^2$. The difficulty is that we
do not know how the log canonical resolution of the ideal sheaf
$\{z\in \mathbb{C}^n:F(z)=F(z_0)\}$ depends on the parameter $z_0$.
On the other hand, the parameter dependence of the complex
singularity exponent is clear from the work of
Demailly-Koll$\acute{\rm a}$r (see also [19] for weaker results).

\medskip

\textbf{Proposition 5$^{[12]}$.}\ {\it Let $M$ be a complex
manifold. Let ${\cal P}(M)$ be the set of locally $L^1$ psh
functions on $M$, equipped with the topology of $L^1$ convergence on
compact subsets. Let $p\in M$ and $\varphi\in {\cal P}(M)$ be given.
If $c<c_p(e^{-\varphi})$ and $\psi$ converges to $\varphi$ in ${\cal
P}(M)$, then $e^{-c\psi}$ converges to $e^{-c\varphi}$ in $L^2$ norm
over some neighborhood $V$ of $p$. }

\medskip

Theorem 1 shows $b_{p_0}(\Omega_F\cap U)=2+c_{z_0}(|F-F(z_0)|)$ for
all $p_0=(z_0,w_0)\in \partial \Omega_F\cap U$, while Proposition 5
implies that for any $c<c_0(|F|)$  there exists a neighborhood $V$
of $0\in \mathbb{C}^{n}$ such that $|F-F(z_0)|^{-c}$ is $L^2$ on $V$
provided $|z_0|$ sufficiently small, consequently,
$c_{z_0}(|F-F(z_0)|)\ge c$ and it follows that the map $p\rightarrow
b_{p}(\Omega_F\cap U)$ is lower semi-continuous on $\partial
\Omega_F\cap U$.

\medskip

\textbf{Question 2.}\ Is the map $p\rightarrow b_p(\Omega)$ lower
semi-continuous on $\partial \Omega$ for any bounded $C^2$
pseudoconvex domain in $\mathbb{C}^n$?

\medskip

\textbf{Remark 8.}\ We do not know whether there exists a bounded
$C^2$ pseudoconvex domain such that the holomorphic sectional
curvature of the Bergman metric is unbounded.

\medskip

\vspace{3mm}\textbf{Acknowledgements}%\dawuhao\.
 \  The authors would like to thank the referee for bringing our
attention that Herbort$^{[20]}$ found recently an example of
pesudoconvex domain whose holomorphic sectional curvature of the
Bergman metric is unbounded.

\end{document}